\documentclass[10pt]{amsart}
\usepackage{amsmath,amsfonts,amssymb,amsxtra,latexsym,amscd,enumerate,amsthm}

\usepackage{latexsym}
\usepackage{epsfig}

\def\t{\theta}

\def\a{\alpha}
\def\q{\frac{1}{2}}

\def\e{\varepsilon}

\def\l{\lambda}
\def\D{\Delta}

\def\les{\lesssim}
\def\m{\mu}

\def\R{\mathbb{R}}
\def\C{\mathbb{C}}
\def\Z{\mathbb{Z}}
\def\N{\mathbb{N}}
\def\S{\mathbb{S}}

\def\beq{\begin{equation}}
\def\eeq{\end{equation}}
\def\eq{\Leftrightarrow}

\def\beq{\begin{equation}}
\def\eeq{\end{equation}}

\newtheorem{t1}{Theorem}
\newtheorem{l1}{Lemma}
\newtheorem{p1}{Proposition}

\begin{document}
\title[Global result for Schr\"odinger Maps]{Global results for Schr\"odinger Maps in dimensions $n \geq 3$}

\author{Ioan Bejenaru}
\address{Department of Mathematics, UCLA, Los Angeles CA 90095-1555}
\email{bejenaru@math.ucla.edu}

\begin{abstract} 

We study the global well-posedness theory for the Schr\"odinger Maps equation. We work in $n+1$ dimensions, for $n \geq 3$, and prove a global well-posedness result for small initial data in $\dot{B}^{\frac{n}{2}}_{2,1}$. 

\end{abstract}

\maketitle

\section{Introduction}

In this paper we continue developing the low regularity theory for the Schr\"odinger map equation, as introduced in \cite{be2}. We briefly go over the derivation of the equation. 

If $(M,g)$ is a Riemannian manifold, then the harmonic maps are smooth maps $\phi:\R^{n} \rightarrow M$ which minimize the Lagrangian:

$$L_g^h=\q \int_{\R^n} |\nabla \phi|_g^2 dx$$

Here and throughout the rest of the paper we choose $M=\S^2$ and we identify $\S^2 \setminus \N$ ($\N$ is the north pole) with the Riemannian surface $(\C, g dz d \bar{z})$ by using the stereographic projection:

$$z \in \C \mapsto (\frac{2 \mbox{Re} z}{1+|z|^2}, \frac{2 \mbox{Im} z}{1+|z|^2}, \frac{1-|z|^2}{1+|z|^2}) \in \S^2$$  

\noindent
where the metric is given by $g(z,\bar{z})=(1+|z|^2)^{-2}$. Our functions will be kept in $L^{\infty}$ and this way avoid the problematic issue that this representation has close to the north pole.  For each $t \in \R$, the energy of the map $z: \R^n \times \R \rightarrow \mathbb{C}$ is defined by:

\beq \label{e1}
L^s_g(z(t))=\q \int_{\R^n} \frac{|\nabla z|^2}{(1+|z|^2)^2} dx
\eeq

The Euler-Lagrange equation of this energy functional is given by:

\beq \label{e2}
\sum_{j=1}^{n} (\frac{\partial}{\partial x_j}- \frac{2 \bar{z}}{1+|z|^2} \frac{\partial z}{\partial x_j}) \frac{\partial z}{\partial x_j}=0
\eeq

The Schr\"odinger Map equation is defined as the evolution equation:

\beq \label{e4}
i\frac{\partial z}{\partial t}= \sum_{j=1}^{n} (\frac{\partial}{\partial x_j}- \frac{2 \bar{z}}{1+|z|^2} \frac{\partial z}{\partial x_j}) \frac{\partial z}{\partial x_j}
\eeq

Written in this form, the Schr\"odinger Map equation is a nonlinear Schr\"odinger equation with derivatives and
had been regarded as too difficult to deal with. As a consequence, the most common approach (at least for low regularity theory)
was to rewrite the equation as the Modified Schr\"odinger Map equation and study the system instead; for reference see \cite{n1}, \cite{n2}, \cite{he}, \cite{ka}, \cite{k1}.

For the case of more regular data, the theory for the Schr\"odinger Map equation is richer and goes beyond the target $\S^2$, covering the more general class of K\"aler manifolds. The main reason for that is that the standard energy methods apply. We refer the reader to the works in \cite{ch}, \cite{dw}, \cite{ks}, \cite{m}.  

We decided to work directly on \eqref{e4} and rewrite it as:

\beq \label{e5}
i z_{t} - \Delta z=  \frac{2 \bar{z}}{1+|z|^2} (\nabla z)^2
\eeq

We have explained thoroughly in \cite{be2} that the derivative part of the nonlinearity for this equation satisfies a null condition. This motivated us to treat \eqref{e5} as it is and not pursue the usual approach via the Modified Schr\"odinger Map equation.

We look at the more general semilinear equation:

\begin{equation} \label{E}
\begin{cases}
\begin{aligned} 
iu_{t}-\Delta u &= Q(u,\bar{u}) (\nabla u)^2, \   t\in \mathbb{R}, x\in \mathbb{R}^{n} \\
u(x,0) &=u_{0}(x)
\end{aligned}
\end{cases}
\end{equation}

\noindent
where $u: \mathbb{R}^{n} \times \mathbb{R} \rightarrow \mathbb{C}$ and $Q: \C \times \C \rightarrow \C$ is analytic. We denote by $N(u)=Q(u,\bar{u}) (\nabla u)^2$ the nonlinearity of our problem. Our nonlinearity in (\ref{e5}) can be written this way as long as we keep $||z||_{L^{\infty}}$ small and we will make sure this is the case.

The scaling for \eqref{E} is $s_c=\frac{n}{2}$. In \cite{be2}, the author established the local well-posedness theory for \eqref{E} for $u_0 \in H^s$ where $s > \frac{n}{2}$. 

Then the author learned that, at about the same time, Ionescu and Kenig have obtained a similar result assuming that $s > \frac{n}{2}+\q$, see \cite{ik2}. The approaches in \cite{ik2} and \cite{be2} are similar up to the point of iterating
the equation, but the spaces in which the iteration is performed are not the same. The author learned from \cite{ik2} the usefulness of 
norms of type $L^{\infty}_{x_\t}L^{2}_{y_\t,t}$, the maximal function type estimate and the related structures in dealing with this problem. These norms proved to be more robust 
in the critical case for a simple reason: they are ready to fit a global in time structure. The wave-packet structure used by the author in 
\cite{be2} was tied up to a local in time approach. 

The $L^{\infty}_{x_\t}L^{2}_{y_\t,t}$ structure has been successfully used by the above authors in dealing with the Benjamin-Ono equation, see \cite{ik1}.

Having these additional structures at his disposal, the author decided to question what happens with the equation \eqref{E} at scaling.

Inspired by the work of Tataru for wave-maps, see \cite{ta}, we decided to investigate the issue of global well-posedness for \eqref{E} for small initial data $u_0 \in \dot{B}^{2,1}_{\frac{n}{2}}$. The advantage is that $\dot{B}^{2,1}_{\frac{n}{2}} \subset L^{\infty}$. Nevertheless  $\dot{B}^{2,1}_{\frac{n}{2}}$ and $\dot{H}^{\frac{n}{2}}$ enjoy the same scaling. As a consequence our problem is a local one and does not see the geometry of the sphere. 

In this setup, Ionescu and Kenig have obtained a positive result for small initial data in $\dot{B}_{\frac{n}{2}}^{2,1}$ for $n \geq 3$, see \cite{ik3}. At the same time we announced a similar result, see Theorem \ref{T} bellow.  Ionescu pointed out a gap in the argument in the preprint we made available. We worked out the missing arguments, plus the necessary adjustments and (about two moths later since) we have the new preprint.
 
From the preprint \cite{ik3}, we learned that Ionescu and Kenig that iterated the equation in the same spaces they have introduced in \cite{ik2}. Hence our approach and the approach of Ionescu and Kenig is similar from this point of view. Nevertheless, there are
significant differences in the way the rest of the argument is carried out in the two papers. 

We end the chronology of the events by mentioning that very recently, Bejenaru, Ionescu and Kenig have obtained a global result for the Schr\"odinger Maps (the setup of the equation is different than the one in this paper) with smooth initial data provided that it has small $\dot{H}^{\frac{n}{2}}$ size, assuming that $n \geq 4$, see \cite{bik}.

\

We briefly describe our (standard) strategy in dealing with the problem. We introduce the inhomogeneous Schr\"odinger equation:

\begin{equation} \label{III}
\begin{cases}
\begin{aligned} 
iu_{t}-\Delta u &= f, \   t\in \mathbb{R}, x\in \mathbb{R}^{n} \\
u(x,0) &=g(x)
\end{aligned}
\end{cases}
\end{equation}

We seek for spaces of functions $Z^{\frac{n}{2}}$ and $W^{\frac{n}{2}}$ with the following properties:

\begin{itemize}

\item{(\it linear property)} \rm The solution $u$ to (\ref{III}) satisfies

\beq \label{le}
||u||_{Z^{\frac{n}{2}}} \les ||g||_{\dot{B}^{2,1}_{\frac{n}{2}}}+||f||_{W^{\frac{n}{2}}}
\eeq

\item{ \it (nonlinear estimate)} \rm  $N$ has the mapping property

$$ N: Z^{\frac{n}{2}} \rightarrow W^{\frac{n}{2}} \ \ \mbox{is Lipschitz continuous}$$

\end{itemize}

Proving the nonlinear estimates involves solving two problems. Since (some of) our spaces are not closed under conjugation we need to define the conjugate of a space $X$, $\bar{X}$, by:

\beq \label{C}
u \in \bar{X} \eq \bar{u} \in X
\eeq

We should prove that

\beq \label{A}
Z^{\frac{n}{2}}+\bar{Z}^{\frac{n}{2}} \ \ \mbox{is an algebra}
\eeq

\noindent
and that $Z^{\frac{n}{2}}+\bar{Z}^{\frac{n}{2}}$ leaves $W^{\frac{n}{2}}$ invariant under multiplication:

\beq \label{M}
(Z^{\frac{n}{2}}+\bar{Z}^{\frac{n}{2}}) W^{\frac{n}{2}} \subset W^{\frac{n}{2}}
\eeq

Finally we need to prove the bilinear estimate. If we denote by $B(u,v)=\nabla{u} \cdot \nabla{v}$ we have to show that:

\beq \label{bee}
||B(u,v)||_{W^{\frac{n}{2}}} \les ||u||_{Z^\frac{n}{2}} ||v||_{Z^{\frac{n}{2}}}
\eeq

Proving the nonlinear estimate amounts to proving (\ref{A})-(\ref{bee}). In order to prove a result for smoother initial data,
we need (smoother) spaces $Z^s$ and $W^s$ with the following properties: 

\begin{itemize}

\item{(\it linear property)} \rm The solution $u$ to (\ref{III}) satisfies

\beq \label{sm0}
||u||_{Z^{s}} \les ||g||_{\dot{H}^{s}}+||f||_{W^s}
\eeq

\item{ \it (nonlinear estimate)} \rm  which is reduced to the estimates:

\end{itemize}

\beq \label{sm1}
||B(u,v)||_{W^s} \les ||u||_{Z^s}|| v||_{Z^{\frac{n}{2}}}+ ||v||_{Z^s}|| u||_{Z^{\frac{n}{2}}}
\eeq

\beq \label{sm2}
||u \cdot v||_{Z^s+\bar{Z}^s} \les ||u||_{Z^s+\bar{Z}^s}|| v||_{Z^{\frac{n}{2}}+\bar{Z}^{\frac{n}{2}}}+ ||v||_{Z^s+\bar{Z}^s}|| u||_{Z^{\frac{n}{2}}+\bar{Z}^{\frac{n}{2}}}
\eeq

\beq \label{sm3}
||u \cdot v||_{W^s} \les ||u||_{Z^{\frac{n}{2}}+\bar{Z}^{\frac{n}{2}}}|| v||_{W^s}+ ||u||_{Z^s+\bar{Z}^s}||v||_{W^\frac{n}{2}}
\eeq

Once we have the linear property and the nonlinear estimate (for the critical regularity and for higher one), then a standard fixed point argument (see for instance \cite{ta}) gives the us the main result:

\begin{t1} \label{T}

a) Assume that $n \geq 3$. Then there exists $\e, \delta > 0$ such that for every $u_{0}$ with $||u_{0}||_{\dot{B}_{\frac{n}{2}}^{2,1}}  < \delta$ , the problem (\ref{E}) has a global solution $u \in C_t \dot{B}_{\frac{n}{2}}^{2,1}$ which is unique in $\{u \in Z^{\frac{n}{2}}: ||u||_{Z^{\frac{n}{2}}} \leq \e\}$. Furthermore, the solution has Lipschtiz dependence on the initial data.

b) If in addition $u_0 \in \dot{H}^s$ with $s > \frac{n}{2}$, then $u \in C_t \dot{H}^s$ and

\[
||u(t)||_{\dot{H}^s} \les ||u_0||_{\dot{H}^s}
\]

\end{t1}

This result is, essentially, the equivalent of the result obtained by Tataru for the wave maps, see \cite{ta}. The result in \cite{ta} is available only for dimension $n \geq 4$ since the geometry of the problem in different. The free solutions for the Schr\"odinger equation live on the paraboloid while those for the wave equation live on the cone. The paraboloid has nonzero curvature, while the cone contains straight lines. This way the Fourier transform of the standard measure on the paraboloid decays faster than its equivalent on the cone. This is what makes our result available in dimension $3$ also. 

Nevertheless the result in dimension $2$ is still open and most likely much harder. The reason is that we lose the $L^{2}_{x}L^{\infty}_{x',t}$ estimate which is crucial for most of the argument.

The rest of the paper is organized as follows. Next section is dedicated to the basic notations and definitions. 
Section \ref{LI} deals with two problems. It derives certain useful norms for solutions for the Schr\"odinger equation and
establishes linear estimates with respect to those norms. Section \ref{DF} defines the spaces we work with as suggested by 
the estimates established in the previous section and, as a consequence, we obtain \eqref{le} for free. Sections \ref{BE} and \ref{AP} deal
with the bilinear estimates, respectively with the algebra type properties, see \eqref{A} and \eqref{M}. In the last section we sketch an argument 
for the case of smoother initial data.  In the Appendix we deal
with an estimate which we decided to postpone in order not to distract the attention of the reader from our main purpose. 

\vspace{.1in}

The author thanks Ionescu and Kenig for making the preprint \cite{ik2} available. The author thanks Ionescu for pointing out a gap in the first draft of this paper and for several other useful commentaries.
The author thanks Daniel Tataru and Terry Tao for useful discussions and encouragement with this project.

\section{Notations and definitions}

Throughout this paper $A \lesssim B$ means $A \leq CB$ for some constant $C$ which is independent of any possible variable in our problem. We say $A \approx B$ if $A \leq CB \leq C^{2}A$ for the same constant $C$.

On the Fourier side, we define the norm of $(\xi,\tau)$ by $|(\xi,\tau)|^{2}=|\tau|+\xi^2$. An important quantity in our paper is $|\tau-\xi^2|$ to which we refer as the modulation. Throughout the paper $\l$ (and the similar entities) and $d$ are dyadic values $2^{j}$ with $j \in \Z$. Using the sum $\sum_\l$ means that we sum over the full range $2^{j}, j \in \Z$, unless other restrictions are specified.

We define the sets:

\[
A_\l=\{ \frac{\l}{2} \leq |(\xi,\tau)| \leq 2 \l \} \ \ \ \mbox{and} \ \ \ B_d=\{ \frac{d}{2} \leq |\tau-\xi^2| \leq 2 d \}
\]

\noindent 
and $A_{\l,d} = A_\l \cap B_d$. Notice that $A_{\l,d}=\emptyset$ if $d \geq 10 \l^2$.

We can construct in a standard way $s \in C_{0}^{\infty}((\q,2))$ such that:

\[
\sum_{\l} s(\frac{x}{\l})=1, \ \forall x \ne 0
\]

We define $s_{\l}(x)=s(\frac{x}{\l})$ and $s_{\leq \l}(x)=\sum_{\l' \leq \l} s_{\l'}(x)$.  Then we introduce the localization operators $S_\l$, $S_{\leq \l}$, $M_{d}$, $M_{\leq d}$ with symbols

\[
s_{\l}(\xi,\tau)=s_\l(|(\xi,\tau)|), \ s_{\leq \l}(\xi,\tau)=s_{\leq \l}(|(\xi,\tau)|) 
\]

\[
m_{d}(\xi,\tau)=s_d(|\tau-\xi^2|), \ m_{\leq d}(\xi,\tau)=s_{\leq d}(|\tau-\xi^2|)
\]

 We define their composition $S_{\l,d}=S_\l M_d$, $S_{\l, \leq d}=S_\l M_{\leq d}$

\[
f_\l=S_\l f \  ,  \  f_{\l,d}=S_{\l,d} f , \ f_{\l, \leq d}=S_{\l, \leq d}f 
\] 

\[
f_{\cdot, d} = M_d f, \ f_{\cdot, \leq d} = M_{ \leq d} f
\]

Occasionally we may need to localize only with respect to the $\xi$ variable, and this is why we introduce also $s_\l(D)$ the 
operator with symbol $s(\frac{|\xi|}{\l})$. 

If $f$ is a distribution on $\R^n$ such that 

\[
f=\sum_{\l} f_{\l}
\]

\noindent
in a distributional sense, we define $\dot{B}^{2,1}_{\frac{n}{2}}$ by 

\[
||f||_{\dot{B}^{2,1}_{\frac{n}{2}}} = \sum_{\l} \l^{\frac{n}{2}} ||f_\l||_{L^{2}}
\]

Here by $f_\l$ we mean $s_\l(D)f$ (we tolerate the fact that above $f_\l$ was defined in the space-time context). As a consequence it follows that $\hat{f} \in L^1$, therefore $f$ has zero limit at $\infty$. 

Throughout the paper by $\hat{f}$ we mean the space-time Fourier transform of $f$, unless otherwise specified. 

For any $b \in \R$ and $1 \leq p \leq \infty$ we define $X_\l^{0,b,p}$ to be the completion of the space of Schwartz functions supported at frequency $\l$ whose norm

\[
||f||^p_{X^{0,b,p}_\l}=\sum_d d^{p b} ||f_{\cdot,d}||^p_{L^{2}} 
\]

\noindent
is finite. From Plancherel we notice the dual relationship:

\beq \label{dual}
sup \{ \langle \phi,\psi \rangle: ||\psi||_{X_\l^{0,\q,1}}=1 \} =||\phi||_{X_\l^{0,-\q,\infty}}
\eeq

\noindent
where $\langle \phi,\psi \rangle$ is the usual inner product:

\[
\langle \phi,\psi \rangle=\int \phi \bar{\psi} dx dt
\]

The spaces $X_{\l}^{0,b,1}$ are atomic spaces whose atoms are functions with Fourier transform supported at frequency $\l$ and modulation $d$ and whose $L^{2}$ norms are of order $O(d^{-b})$. As a consequence, functions in $X_\l^{0,\q,1}$ are uniquely defined modulo solutions of the homogeneous Schr\"odinger equation. We can remedy this if we include the $L^2$ solutions of the homogeneous Schr\"odinger equation in the class of acceptable atoms in $X_\l^{0,\q,1}$ with their $L^2$ norm. 

For a more detailed description of the $X^{s,b,p}$ spaces we refer the reader to \cite{ta} and \cite{tao}.

Besides Bourgain type spaces, we involve space-time structures. Even though these spaces will be defined later in the paper, it is more appropriate
to prepare some of the ingredients here.  

For each $\t=(\t_1,...,\t_{n-1}) \in [0,2\pi] \times [0,\pi]^{n-2}$ we associate the vector $v_{\t}=(1,\t)=(1,\t_1,...,\t_{n-1}) \in S^{n-1}$; the representation of $v_\t$ is meant in polar coordinates. The set $\{\xi \in \R^{n}: \xi \cdot v_{\t} \geq \frac{|\xi|^2}{2} \}$ which collects the vectors "pointing" in the same direction as $v_\t$ will play an important role in our argument.

We would like to cover $\R^{n}$ with a minimal number of such sets. Since $\xi \cdot v_\t \geq \frac{|\xi|^2}{2}$ is equivalent to $\angle{(\xi,v_{\t})} \leq \frac{\pi}{3}$, it follows that we can select a set $\Theta \subset [0,2\pi] \times [0,\pi]^{n-2}$ of cardinality at most $10^{n}$ such that

 \[
 \R^{n}=\cup_{\t \in \Theta} A_{\t}
 \]

\noindent 
and $\xi \cdot v_\t \geq \frac{|\xi|^2}{2}$ for all $\xi \in A_\t$. We can also impose a condition on $\Theta$ so that we avoid duplicates of $A_{\t}$. For instance we can ask $\angle{v_\t,v_{\t'}} \geq \frac{\pi}{10}$ for any $\t \ne \t'$. 

 Whenever we refer to $\t$ we implicitly understand that $\t \in \Theta$. For each $\t$ we associate the orthogonal coordinates on the physical side $(x_{\t},x_{\t}')$ where:

\[
x_{\t}=x \cdot v_{\t}
\]   

We denote by $(\xi_{\t},\xi'_{\t})$ the corresponding Fourier variables. The symbol $\tau-\xi^{2}$ does not change under this change of coordinates and that $d\xi=d\xi_{\t} d\xi_{\t}'$, therefore norms on the Fourier side do not change either.  

We need an operator $S_\t$ to localize, on the frequency side, the angular variable in $A_\t$. We can build a system of functions $h_\t$ such that $h_\t$ is a smooth approximation of $\chi_{A_\t}$ and $(h_\t)_{\t \in \Theta}$ form a partition of unity. Then we define $S_\t$ the operator with symbol $h_\t(\xi)$.     

For each $\t$ we introduce two additional types of sets:

\[
B_{\t}=\{(\xi,\tau) \in R^{n+1}: \xi \in A_{\t}, |\tau - \xi^{2}| \leq \frac{\xi^{2}}{10}\}
\]

\[
\tilde{A}_{\t}=\{(\xi_{\t}',\tau): (\xi,\tau) \in B_{\t} \}
\]

$B_{\t}$ proves to be the useful lift up  of $A_{\t}$ in $\R^{n+1}$. 
$A_\t$ is the projection of $B_{\t}$ onto the plane $\tau=0$ 
(along the $\tau$ direction) and $\tilde{A}_\t$ is the projection of 
$B_{\t}$ on the plane $\xi_\t=0$ (along the vector $v_{\t}$). 

The important facts about $(\xi,\tau) \in B_{\t}$ are:

\[
\tau \geq \xi^{' 2}_{\t}, \  \ \xi_{\t} \approx |\xi| \ \ \mbox{and} \ \ \xi_{\t}+\sqrt{\tau-\xi_{\t}^{' 2}} \approx |\xi|
\]

\section{Linear Estimates} \label{LI}

We want to investigate the linear Schr\"odinger equation:

\begin{equation} \label{II}
\begin{cases}
\begin{aligned} 
iu_{t}-\Delta u &= f, \   t\in \mathbb{R}, x\in \mathbb{R}^{n} \\
u(x,0) &=g
\end{aligned}
\end{cases}
\end{equation}

For this equation we have the Strichartz estimates:

\beq \label{ST}
||u||_{L^{p}_{t}L^{q}_{x}} \les ||f||_{L^{\tilde{p}}_{t}L^{\tilde{q}}_{x}}+||g||_{L^2_x}
\eeq

\noindent
where $(p,q)$ is any admissible pair, i.e. a pair satisfying $\frac{n}{q}+\frac{2}{p}=\frac{n}{2}$, $2 < p \leq \infty$, and $(\tilde{p},\tilde{q})$ is any dual pair to an admissible pair, i.e. there is an admissible pair $(\tilde{p},\tilde{q})$ such that $\frac{1}{\tilde{p}}+\frac{1}{\tilde{p}}=\frac{1}{\tilde{q}}+\frac{1}{\tilde{q}}=1$. For reference, see \cite{gv}.

These estimates are derived by seeing $t$ as the evolution direction of the equation. We can obtain useful information about the solution of \eqref{II} if we change the evolution direction of the equation. 

The solution for the free evolution ($f=0$) is $u=e^{-it\D}g$ and $\hat{u}$ lives on the paraboloid $\tau=\xi^2$ since
 
 \[
 (\tau-\xi^2)\hat{u}=0
 \]

We decompose the symbol $\tau-\xi^2$ in 

\beq \label{decs}
\tau-\xi^2=\tau-\xi_{\t}^{' 2}-\xi_{\t}^2=(\sqrt{\tau-\xi_{\t}^{' 2}}-\xi_{\t})(\sqrt{\tau-\xi_{\t}^{' 2}}+\xi_{\t}) 
\eeq

\noindent
and notice that if we localize $\hat{g}$ in $A_{\t}$ it follows that $\hat{u}$ is localized in $B_{\t}$ where we have $(\sqrt{\tau-\xi_{\t}^{' 2}}+\xi_{\t}) \approx |\xi| > 0$. Therefore our equation is equivalent to

\[
(\sqrt{\tau-\xi_{\t}^{' 2}}-\xi_{\t})\hat{u}=0
\]

\noindent
which can be seen now as an evolution equation in the direction of $x_{\t}$:

\beq \label{hh}
(\partial_{x_{\t}}+i\sqrt{\tau-\xi^{' 2}_{\t}})\hat{u}=0
\eeq 

\noindent
where $\hat{u}$ is now the Fourier transform of $u$ with respect to $x'_{\t},t$ variables. We define $a(D_{x'_{\t}},D_t)$ the pseudo-differential operator with symbol:

\beq \label{defa}
a(\xi'_{\t},\tau)=\sqrt{\tau-\xi_{\t}^{' 2}}
\eeq

The solution \eqref{hh} is $u(x_{\t},x',t)=e^{-ix_{\t} a(D_{x_{\t}'},D_t)} \psi$ where $\psi=u(0,x_{\t}',t)$ is the initial data of \eqref{defa}. One requires a bit of care when transferring the initial data from \eqref{II} to \eqref{hh}, but this will not be an issue. In any case one expects control on the norm $||u||_{L^{\infty}_{x_{\t}}L^{2}_{x_{\t}',t}}$.

The above change of the evolution coordinate is another way to say that we can write the equation of the paraboloid  in two different ways:  $\tau=\xi^2$ and $\xi_\t=\sqrt{\tau-\xi'^2_\t}$ (as long as we stay localized in $B_{\t}$). 

Now we turn our attention to the inhomogeneous problem which, after taking a complete Fourier transform, becomes:

\[
 (\tau-\xi^2)\hat{u}=\hat{f}
 \]

Using the decomposition  in \eqref{decs},we rewrite our equation as:

\[
(\sqrt{\tau-\xi_{\t}^{' 2}}-\xi_{\t})\hat{u}=(\sqrt{\tau-\xi_{\t}^{' 2}}+\xi_{\t})^{-1} \hat{f}
\]

\noindent

or

\beq \label{ih3}
(\partial_{x_{\t}}+i\sqrt{\tau-\xi'^2})\hat{u}=\int e^{ix_{\t} \xi_{\t}} (\sqrt{\tau-\xi'^2}+\xi_{\t})^{-1} \hat{f} d\xi_{\t}
\eeq

\noindent
where, in the last equation, the Fourier transform of $u$, $\hat{u}$, is taken with respect to the variables $x_{\t}',t$.  This is an evolution equation in the direction of $x_{\t}$ and a reasonable estimate to expect is that $||f||_{L^{1}_{x_{\t}}L^{2}_{x_{\t}',t}}$ gives control on $||u||_{L^{\infty}_{x_{\t}}L^{2}_{x_{\t}',t}}$. If  $\hat{f}$ is localized in $B_{\t}$ and at frequency $\l$, then $(\sqrt{\tau-\xi_{\t}^{' 2}}+\xi_{\t})^{-1} \approx \l^{-1}$, hence, recalling \eqref{ih3}, we actually expect $||u||_{L^{\infty}_{x_{\t}}L^{2}_{x_{\t}',t}}$ to be controlled by $\l^{-1}||f||_{L^{1}_{x_{\t}}L^{2}_{x_{\t}',t}}$. In other words, with respect to these norms, we record a global smoothing effect: our equation recovers a full derivative from the inhomogeneity and this property is known to be crucial in dealing with nonlinearities with derivatives.

The above calculus was developed by the author in \cite{be1}.

In what follows we make all the above heuristics rigorous. In dealing with estimates it is preferable to make use of the scaling properties of the equation. If $u$ is a solution of \eqref{I} then $u_{\l}= u(\frac{x}{\l},\frac{t}{\l^2})$ is a solution for the same equation where $f_{\l}= \l^{-2} f(\frac{x}{\l},\frac{t}{\l^2})$ and $g_{\l}=g(\frac{x}{\l})$. This allows us to rescale the estimates localized in frequency and prove them at frequency of size $\approx 1$. 

Our first result deals with the estimates for the free evolution.

\begin{p1}
a) If $u_0 \in L^{2}_x$ localized at frequency $\l$ and in $A_{\t}$ then:

\beq \label{be}
||e^{-it\D}u_0||_{L^{\infty}_{x_{\t}}L^{2}_{t,x_{\t}'}} \les \l^{-\frac{1}{2}} ||u_0||_{L^2}  
\eeq

b) If $u_0 \in L^{2}_{x_{\t}',t}$ localized at frequency $\l$ and in $\tilde{A}_{\t}$ then:

\beq \label{be2}
||e^{-ix_{\t} a(D_{x_{\t}'},D_t)}u_0||_{L^{\infty}_{x_{\t}}L^{2}_{t,x_{\t}'}} \les ||u_0||_{L^2_{x_{\t}',t}}  
\eeq

\beq \label{be1}
||e^{-ix_{\t} a(D_{x_{\t}'},D_t)}u_0||_{L^{\infty}_{t}L^{2}_x} \les \l^{\q} ||u_0||_{L^2_{x_{\t}',t}}  
\eeq

\beq \label{be3}
||e^{-ix_{\t} a(D_{x_{\t}'},D_t)}u_0||_{L^{\frac{2(n+2)}{n}}_{t,x}} \les \l^{\q} ||u_0||_{L^2_{x_{\t}',t}}  
\eeq

In \eqref{be1} one gets the same control over the norm in $C_t L^2_x$. 

\end{p1}

\begin{proof}

By rescaling we can assume that $\l \approx 1$. The problem is equivalent to showing:

$$||\int e^{i x \cdot \xi} e^{i t \xi^2} \psi(\xi) d\xi||_{L^{\infty}_{x_{\t}}L^{2}_{t,x_{\t}'}} \les ||\psi||_{L^2}$$

\noindent
where $\psi$ (which replaces $\hat{u}_0$) is localized at $|\xi| \approx 1$ and in $A_{\t}$. We perform the change of variables $\xi_{\t}=\sqrt{\tau-\xi_{\t}^{' 2}}$ and $\xi'_{\t}=\xi'_{\t}$ in the left term above and then estimate

\[
\begin{split}
&||\int e^{i (x_{\t}' \cdot \xi'_\t+ t\tau)} e^{i x_{\t} \sqrt{\tau-\xi_{\t}^{' 2}} } \psi(\sqrt{\tau-\xi_{\t}^{' 2}},\xi'_{\t}) (2\sqrt{\tau-\xi'^{2}})^{-1} d\xi'_{\t} d\tau||_{L^{\infty}_{x_{\t}}L^{2}_{t,x_{\t}'}} \\
&\les || e^{i x_{\t} \sqrt{\tau-\xi^{' 2}_{\t}} } \psi(\sqrt{\tau-\xi_{\t}^{' 2}},\xi'_{\t}) (2\sqrt{\tau-\xi_{\t}^{' 2}})^{-1}||_{L^{2}_{\xi'_{\t},\tau}} \\
&\approx || \psi(\sqrt{\tau-\xi^{' 2}_{\t}},\xi'_{\t}) (2\sqrt{\tau-\xi^{' 2}_{\t}})^{-1}||_{L^{2}_{\xi'_{\t},\tau}} \approx |||\xi_{\t}|^{-\q} \psi(\xi_{\t},\xi'_{\t})||_{L^{2}_{\xi}} \approx ||\psi||_{L^{2}}
\end{split}
\]

In the last line we have used the localization property of $\psi$ via the fact that $|\xi_{\t}| \approx 1$ for $\xi \in A_{\t}$. 

\eqref{be2} is obvious, while \eqref{be1} can be obtained by going reverse way in the previous estimate. After rescaling, \eqref{be3} is the standard Strichartz estimate, since at frequency $1$ and in $\tilde{A}_\t$ the induced measures on $\tau=\xi^2$ and $\xi_{\t}=\sqrt{\tau-\xi_{\t}^{' 2}}$ are the same.

\end{proof}

Next, we derive the maximal function estimate for the free solutions. This estimate plays a crucial role in the bilinear estimates and as we can see bellow we do not have the estimate in dimension two. 

\begin{p1}
Assume $u_0 \in L^{2}_x$ localized at frequency $\l$. Then for every $\t$:

\beq \label{se}
||e^{-it\D}u_0||_{L^{2}_{x_{\t}}L^{\infty}_{t,x_{\t}'}} \les \l^{\frac{n-1}{2}} ||u_0||_{L^2}  
\eeq

If $u_0 \in L^{2}_{x_{\t'}',t}$ is localized at frequency $\l$ and in $\tilde{A}_{\t'}$ then for every $\t$:

\beq \label{se2}
||e^{-ix_{\t} a(D_{x_{\t}'},D_t) }u_0||_{L^{2}_{x_{\t}}L^{\infty}_{t,x_{\t}'}} \les \l^{\frac{n}{2}} ||u_0||_{L^2_{x_{\t'}',t}}  
\eeq

\end{p1}

\bf Remark. \rm The first estimate gives us control over the norm $||e^{-it\D}u_0||_{L^{2}_{x_{\t}}L^{\infty}_{t,x_{\t}'}}$ regardless of the choice of direction given by $\t$. The second estimate gives us a similar control even if we have control on the input in a preferential direction, the one given by $\t'$. 

\begin{proof}
By rescaling we can assume that $\l \approx 1$. Then \eqref{se} is equivalent to the estimate:

\beq \label{int1}
||\int e^{i x \cdot \xi} e^{i t \xi^2} \psi(\xi) d\xi||_{L^{2}_{x_{\t}}L^{\infty}_{t,x_{\t}'}} \les ||\psi||_{L^2}
\eeq

Here $\psi=\hat{u}_0$. By the $TT^{*}$ argument this is equivalent to showing that:

$$||\int e^{i x \cdot \xi} e^{i t \xi^2} s(|\xi|) d\xi||_{L^{1}_{x_{\t}}L^{\infty}_{t,x_{\t}'}} \les 1$$

Recall that $s(|\xi|)$ localizes at $|\xi| \approx 1$. If $d \mu$ is the measure on the paraboloid $\tau=\xi^2$ where $|\xi| \approx 1$, the above estimate is equivalent to:

\beq \label{se1}
||\widehat{d \mu}(x,t)||_{L^{1}_{x_{\t}}L^{\infty}_{t,x_{\t}'}} \les 1
\eeq

 From the standard theory of Fourier transforms of measures supported on surfaces, see \cite{st}, we have the estimate:

$$|\widehat{d \mu}(x,t)| \les \langle (x,t)\rangle^{-\frac{n}{2}} $$

Therefore as long as $n \geq 3$ this decay is enough to guarantee the estimate \eqref{int1} for any choice of $x=(x_{\t},x_{\t}')$. 

For \eqref{se2} we rescale to bring the estimate to $\l \approx 1$. Then we notice that we deal with the same problem as in \eqref{se} since the surface $\xi_{\t}=\sqrt{\tau-\xi_{\t}^{' 2}}$ is the same with $\tau=\xi^2_{\t}+\xi^{' 2}_{\t}$ and the induced measures are the equivalent (since we work at $|\xi| \approx 1$).   

\end{proof}

Now that we have the estimates for the free evolution, the next natural step is to derive the corresponding ones for the inhomogeneous equation. There is one additional technicality involved when dealing with the inhomogeneous equation. In addition to the estimates for the solution of 

\begin{equation} \label{I}
\begin{cases}
\begin{aligned} 
iu_{t}-\Delta u &= f, \   t\in \mathbb{R}, x\in \mathbb{R}^{n} \\
u(x,0) &=0
\end{aligned}
\end{cases}
\end{equation}

\noindent
we need also estimates for the truncated solutions with respect to the modulation. 
The following Lemma provides one of the ingredients for this process.

\begin{l1} \label{cons} 

a)The operators $S_{\l,\leq d}$ are bounded on $L^{p}_t L^{2}_x$, for $1 \leq p \leq \infty$.  

b)The operators $S_{\l,\leq d} S_{\t}$ are bounded on $L^{p}_{x_\t} L^{2}_{x_\t',t}$, for $1 \leq p \leq \infty$. 

\end{l1}

\begin{proof}

The argument here is similar to the one in \cite{ta}. In fact, part a) has a complete proof there, the only adjustment needed being 
the replacement of the cone with the paraboloid.

We prove part b) instead, since this is requires a bit of care. One can easily adapt the argument below for part a) too.
Assume first that $d \leq \frac{\l^2}{20}$, so that things are localized, on the Fourier side, in $B_\t$. We denote by $\tilde{s}=s_{\l,d} \cdot s_\t$. For fixed $\xi_{\t}'$ and $\tau$, $\tilde{s}$ is supported in the region $|\xi_\t-\sqrt{\tau-\xi'_\t}| \les \l^{-1} d$ and satisfies:

\[
|D^{\a}_{\xi_\t} \tilde{s}| \les c_\a |\l^{-1} d|^{-|\a|}
\]

Then, the inverse Fourier transform of $\tilde{s}$ with respect to $x_\t$, $K(x_\t,\xi_\t',\tau)$ satisfies:

\[
K(x_\t,\xi_\t',t)=e^{ix_\t \sqrt{\tau-\xi_\t^{'2}}} L(x_\t,\xi_\t',\tau)
\]

where

\[
|L(x_\t,\xi_\t',t)| \leq c_{N} \frac{\l^{-1} d}{(1+ |x_\t|\l^{-1} d)^{N}}
\]

Hence $K(x_\t,\xi_\t',t) \in L^{1}_{x_\t}L^{\infty}_{\xi_\t',\tau}$ and this justifies our claim. If $d \geq \frac{\l^2}{20}$, we split $S_{\l,d} S_\t=S_{\l, \leq \frac{\l^2}{20}} S_{\l,d} S_\t + (1-S_{\l, \leq \frac{\l^2}{20}})S_{\l,d}S_\t $. For $S_{\l, \leq \frac{\l^2}{20}} S_{\l,d} S_\t$ we can apply the previous argument, while for $(1-S_{\l, \leq \frac{\l^2}{20}})S_{\l,d}S_\t $ we can run a simpler argument since the paraboloid plays no role anymore due to the localization. We leave the rest of the details to the reader.

\end{proof}

The result in the above Lemma allows us to conserve $L^{p}_{x_\t} L^{2}_{x_\t',t}$ norms under modulation localizations, as long as we control the angular localization. 

\begin{p1} \label{repeat}

Assume $f \in L^1_{x_{\t}}L^{2}_{t,x_{\t}'}$ and $f$ localized at frequency $\l$ and in $B_\t$. If $u$ is a solution of \eqref{I} then it satisfies 

\beq \label{ie}
||u||_{L^{\infty}_{x_{\t}}L^{2}_{t,x_{\t}'}} \les \l^{-1} ||f||_{L^{1}_{x_{\t}}L^{2}_{t,x_{\t}'}} 
\eeq

\beq \label{ie2}
||u||_{L^{2}_{x_{\t'}}L^{\infty}_{t,\tilde{x}'_{\t'}}} \les \l^{\frac{n}{2}-1} ||f||_{L^{1}_{x_{\t}}L^{2}_{t,x_{\t}'}} 
\eeq

\beq \label{ie3}
||u||_{L^{\infty}_{t}L^{2}_{x}} \les \l^{-\q} ||f||_{L^{1}_{x_{\t}}L^{2}_{t,x_{\t}'}} 
\eeq

\beq \label{ie4}
||u||_{L^{\frac{2(n+2)}{n}}_{x,t}} \les \l^{-\q} ||f||_{L^{1}_{x_{\t}}L^{2}_{t,x_{\t}'}}
\eeq

\eqref{ie2} is valid for any choice of $\t'$. In \eqref{ie3} one gets the same control over the norm in $C_t L^2_x$. 

In addition one obtains all the above estimates for $u_{\cdot, \leq d}$ for any $d$.

\end{p1}

\begin{proof}
We start with \eqref{ie}, then rescale and reduce the problem to the case $\l \approx 1$. Then we recall the derivation of \eqref{ih3} and rewrite our equation as:

\[
(\partial_{x_{\t}}+i\sqrt{\tau-\xi^{'^2}_{\t}})\hat{u}=\int e^{ix_{\t} \xi_{\t}} (\sqrt{\tau-\xi^{' 2}_{\t}}+\xi_{\t})^{-1} \hat{f} d\xi_{\t}
\]

In the above formula, the Fourier transform is taken with respect to the variables $x_\t',t$. We consider a solution of the equation to be:

\[ \label{ps}
u(x_{\t},\cdot)=\int_{-\infty}^{x_{\t}} e^{-i(x_{\t}-s)\sqrt{\tau-\xi'^2}} \int e^{i s \xi_{\t}} (\sqrt{\tau-\xi'^2}+\xi_{\t})^{-1} \hat{f} d\xi_{\t} ds
\]

\noindent
and be aware that we lose track of the initial data.

Since the localization in frequency is in $B_{\t}$ and at frequency $1$, a similar argument to the one in Lemma \ref{cons} shows that $|\mathcal{F}^{-1}((\sqrt{\tau-\xi^{' 2}_{\t}}+\xi_{\t})^{-1} f )|_{L^{1}_{x_\t}L^{2}_{x'_\t,t}} \les ||f||_{L^{1}_{x_\t}L^{2}_{x'_\t,t}}$, therefore we can just ignore the term and work as if:

\beq \label{repr}
u(x_{\t},\xi',\tau)=\int_{-\infty}^{x_{\t}} e^{-i(x_{\t}-s)\sqrt{\tau-\xi^{' 2}_{\t}}}  \hat{f}(s,\xi'_{\t},\tau) ds
\eeq

Using \eqref{be2} gives us the estimate for the non-retarded operator:

\[
||\int_{-\infty}^{\infty} e^{-i(x_{\t}-s)\sqrt{\tau-\xi^{' 2}_{\t}}}  \hat{f}(s,\xi'_{\t},\tau) ds||_{L^{\infty}_{x_{\t}}L^{2}_{t,x_{\t}'}} \les ||f||_{L^{1}_{x_{\t}}L^{2}_{t,x_{\t}'}}
\]

Then we use the Christ-Kiselev Lemma, see \cite{ck}, to claim the retarded estimate \eqref{ie}.

In some sense, one should obtain on  the same line of arguments  \eqref{ie2} (by using \eqref{se2}),  \eqref{ie3} (by using \eqref{be1}) and \eqref{ie4}. Unfortunately we lose track of the retarded variable and this why a direct application of the Christ-Kiselev Lemma cannot be used. Therefore we derive those estimates in the old fashion way via direct estimates.  

We can rewrite  \eqref{repr} in the following form:

\[
u(x_{\t},x'_\t,t)=\int_{-\infty}^{x_{\t}} \int e^{-i[(x_{\t}-s)\sqrt{\tau-\xi^{' 2}_{\t}}+ x'_\t \xi_\t'+ t\tau ]}  \hat{f}(s,\xi'_{\t},\tau) d\xi_\t' d\tau ds
\]

\noindent
from which we conclude with

\[
|u(x_{\t},x'_\t,t)| \les \int_{-\infty}^{+\infty} |\int e^{-i[(x_{\t}-s)\sqrt{\tau-\xi^{' 2}_{\t}}+ x'_\t \xi_\t'+ t\tau ]}  \hat{f}(s,\xi'_{\t},\tau)d\xi_\t' d\tau| ds
\]

\noindent
and this allows us to treat the estimate as a non retarded one and claim \eqref{ie2}, \eqref{ie3} and \eqref{ie4}. 

Note that all the above estimates were derived for a solution of \eqref{I} which did not satisfy the initial data condition. \eqref{ie3} tells us that this solution corresponds to an initial data $u(x,0)$ satisfying an estimate:

\beq
||u(x,0)||_{L^{2}_{x}} \les \l^{-\q} ||f||_{L^{1}_{x_{\t}}L^{2}_{t,x_{\t}'}} 
\eeq

Therefore the true solution of \eqref{I} is $\tilde{u}(x,t)=u(x,t)-e^{-it\D} u(x,0)$. By combining the estimates we have just obtained with the ones for the homogeneous solution $e^{-it\D} u(x,0)$ we obtain the claims for the right solution.

In order to show that besides the $L^{\infty}_t L^2_x$ estimate one obtains a $C_t L^2_x$ estimate, one notices that if  $\hat{f}$  were smooth, then, since it already has compact support, $f$ would be smooth and decaying at infinity, hence one would obtain for free the $C_t L^{2}_x$ estimate for $u$. For general $f$ we can easily construct a sequence $f_n$ such that $\hat{f}_n$ is smooth and supported at frequency $\l$ and $f_n \rightarrow f$ in $L^{1}_{x_\t} L^{2}_{x'_\t,t}$. As a consequence $u_n \rightarrow u$ in $L^{\infty}_t L^2_x$ and since $u_n \in C_t L^2_x$ we obtain that $u \in C_t L^2_x$. 

In the end we need to justify why $u_{\cdot, \leq d}$ satisfies all the estimates. One does not expect 
conservation of all the norms we listed for $u$ under modulation cut-offs. This comes from free from the observation
that $u_{\cdot, \leq d}$ solves the equation:

\[
(i\partial_t - \D)u_{\cdot, \leq d}=f_{\cdot, \leq d}
\]

The result in Lemma \ref{cons} tells us that $||f_{\cdot, \leq d}||_{L^{1}_{x_{1}}L^{2}_{t,x_{\t}'}} \les ||f||_{L^{1}_{x_{1}}L^{2}_{t,x_{\t}'}}$.  Lemma \ref{cons} also tells us that $||u_{\cdot, \leq d}||_{L^{\infty}_t L^{2}_x} \les ||u||_{L^{\infty}_t L^2_x}$ and this gives us an $L^2$ estimate for $u_{\cdot, \leq d}$ at time zero, therefore we can redo the argument above and claim all the estimates for $u_{\cdot, \leq d}$.

\end{proof}

In the next result we would like to claim similar estimates to the one in Proposition \ref{repeat} if we assume that $f \in L^{\frac{2(n+2)}{n+4}}$ instead of $f \in L^1_{x_{\t}}L^{2}_{t,x_{\t}'}$. There is one difficulty we encounter. In Proposition \ref{repeat} we proved that all the
estimates hold true for $u_{\cdot, \leq d}$ by using the fact that $L^1_{x_{\t}}L^{2}_{t,x_{\t}'}$ is stable under modulation cut-offs, see Lemma \ref{cons}. We do not have a similar result for $ L^{\frac{2(n+2)}{n+4}}$.

The way to get around this difficulty is to provide an abstract result saying that even if the right hand side of \eqref{I}
does not carry a norm which is stable under modulation cut-offs, the solution of \eqref{I} does have this property. 
Hence we first list the main result, the version of Proposition \ref{repeat} when $f \in L^{\frac{2(n+2)}{n+4}}$ and then
show, via an abstract result, that we can derive the estimates for $u_{\cdot, \leq d}$ too.

\begin{p1} \label{repeat1}
All the results in Proposition \ref{repeat} hold true if $f \in L^1_{x_{\t}}L^{2}_{t,x_{\t}'}$ is replaced by $f \in L^{\frac{2(n+2)}{n+4}}$ and a factor of $\l^{-\q}$ is taken off from the right hand side of all estimates. 
\end{p1}

\begin{proof} The equivalent of \eqref{ie3} and \eqref{ie4} are the Strichartz estimates, see \eqref{ST}. 
To prove the equivalent of \eqref{ie} we rewrite \eqref{repr} as 

\[
\begin{split}
u(x_{\t},\xi',\tau)&=\int_{-\infty}^{x_{\t}} e^{-i(x_{\t}-s)a(D_{x_{\t}'},D_t)}  f(s,x_\t',t) ds \\
&=e^{-ix_{\t}a(D_{x_{\t}'},D_t)} \int_{-\infty}^{x_{\t}} e^{i s a(D_{x_{\t}'},D_t)}  f(s,x_\t',t) ds
\end{split}
\]

For fixed $x_\t$ we have that:

\[
||\int_{-\infty}^{x_{\t}} e^{i s a(D_{x_{\t}'},D_t)}  f(s,x_\t',t) ds||_{L^{2}} \les ||f||_{L^{\frac{2(n+2)}{n+4}}}
\]

\noindent
this being the dual Strichartz estimate. Then we apply \eqref{be2} to claim the result. 

The proof of \eqref{ie2} is the nontrivial part since we do not have input in a $L^{1}L^{2}$ or output in a $L^{\infty}L^{2}$, these
being the cases which can be treated in a simple manner. We dedicate the Appendix to this estimate.

\end{proof}
 
Below we develop the abstract machinery to allow us to claim the estimates for $u_{\cdot,\leq d}$ in Proposition \ref{repeat1}. The Schr\"odinger equation:

\[
(i\partial_t-\D)u=f
\]

 can be rewritten as

\[
\partial_t (e^{it\D}u)=e^{it\D}f
\]

Then, the solution solution is:

\[
u(t)=\int_{-\infty}^t e^{-i(t-s)\D} f(s) ds=\int_{-\infty}^{\infty} 1_{s \leq t} e^{-i(t-s)\D} f(s) ds
\]

This solution corresponds to a zero initial data at $-\infty$; in the end we explain how we fit a given initial data
(at time $0$) in this argument.

Recalling that $\hat{u}_{\cdot,\leq d}(\xi,\tau)=s_{\leq d}(|\tau-\xi^2|) \hat{u}(\xi,\tau)$ , we have $e^{it\D} u_{\cdot,\leq d}=\phi_{\leq d} * e^{it\D}u$, where $\phi_d=\check{s}_{\leq d}$ is the inverse Fourier transform of $s_{\leq d}$ and the convolution is performed with respect to the $t$ variable.
Therefore we can continue with:

\[
e^{it\D} u_{\cdot, \leq d}(t)= \int \phi_{d}(s) \int_{-\infty}^{\infty} 1_{s' \leq t-s} e^{is'\D} f(s') ds'ds
\]

We can rewrite this as

\[
u_{\cdot, \leq d}(t) = \int \phi_{d}(s) \int_{-\infty}^{\infty} 1_{s' \leq t-s} e^{-i(t-s')\D} f(s') ds'ds =\int \phi_{d}(s) h_s(t) ds
\]

where

\beq \label{dr}
h_s(t)=\int_{-\infty}^{\infty} 1_{s' \leq t-s} e^{-i(t-s')\D} f(s') ds'
\eeq

The operator giving the formula \eqref{dr} is a delayed version of a retarded operator as explained before. The delay is meant in the sense that $1_{s' \leq t}$ is replaced by $1_{s' \leq t-s}$. A closer look at the estimates we provided before shows that they can be easily adapted for:

\[
||h_s||_{X} \les C_{\l} ||f||_{L^{\frac{2(n+2)}{n+4}}}
\]

\noindent
with the implicit constant (see $\les$) independent on $s$. Here $X$ is any of the spaces $L^{2}_{x_\t}L^{\infty}_{x_\t',t}$, $L^{\infty}_{x_\t} L^{2}_{x'_\t,t}$, $L^{\infty}_t L^2_x$, $L^{\frac{2(n+2)}{n}}_{x,t}$. Therefore we obtain that:

\[
||u_{\cdot,\leq d}||_{X} \les C_\l \int |\phi_{d}(s)| ||f||_Y ds = C_\l ||\phi_{d}||_{L^1} ||f||_Y \les C_\l ||f||_{L^{\frac{2(n+2)}{n+4}}}
\]

Since in this process we obtain the $L^{\infty}_t L^2_x$ bound, we can always correct things to fit a given initial data in $L^2_x$. 

This completes missing in the proof of Proposition \ref{repeat1}.

\vspace{.1in}

In the end we describe the connection with the $X^{s,\q,p}_\l$ type spaces. For a function $f \in X_\l^{0,\q,1}$ we use its atomic decomposition:

\[
f=\sum_d f_{\cdot,d} + e^{it\D} \tilde{f}
\]

\noindent
where $f_{\cdot,d} \in X^{0,\q}_{\l,d}$ and $\tilde{f} \in L^{2}$. For each $d$ we have:
\[
f_{\cdot, d}(x,t)=\int e^{i(x \xi + t \tau)} \hat{f}_{\cdot,d}(\xi,\tau) d\xi d\tau= \int e^{i(x \xi + t(\xi^{2}+s))}\hat{f}_{\cdot,d}(\xi,\xi^2+s) d\xi d s=\int e^{its} h_s(x,t) ds
\]

\noindent
where $h_s(x,t)= \int e^{i(x \xi + t\xi^{2})}\hat{f}_{\cdot,d}(\xi,\xi^2+s) d\xi$ is a solution of the homogeneous Schr\"odinger equation with initial data $h_s(x,0)= \int e^{i x \xi}\hat{f}_{\cdot,d}(\xi,\xi^2+s) d\xi$. The range of integration with respect to $s$ is included in  $[4^{-1}d,4d]$, hence:

\[
\int ||h_s(x,0)||_{L^{2}_x} ds \les \int ||\hat{f}_{\cdot,d}(\xi,\xi^2+s)||_{L^{2}_\xi} ds \les d^{-\q} ||\hat{f}_{\cdot,d}(\xi,\xi^2+s)||_{L^{2}_{\xi,s}} \approx ||f_{\cdot,d}||_{X_{\l,d}^{0,\q}}
\]

This allows us to look at $X^{0,\q,1}_\l$ as a $L^1$ superposition of "almost" free solutions and claim all the estimates from before: 

- If $f \in X_{\l}^{0,\q,1}$ and $\hat{f}(\tau, \cdot)$ is localized in $A_{\t}$, for any $\tau$, then:

\beq \label{xsb1es1}
||f||_{L^{\infty}_{x_{\t}}L^{2}_{t,x_{\t}'}} \les \l^{-\q} ||f||_{X_{\l}^{0,\q,1}}
\eeq

-  If $f \in X_{\l}^{0,\q,1}$ then for any $\t$:

\beq \label{xsb1es2}
||f||_{L^{2}_{x_{\t}}L^{\infty}_{t,x_{\t}'}} \les \l^{\frac{n-1}{2}} ||f||_{X_{\l}^{0,\q,1}}
\eeq

\beq \label{xsbstr}
||f||_{L^{p}_{t}L^{q}_{x}} \les ||f||_{X_{\l}^{0,\q,1}}
\eeq

\noindent
for every $(p,q)$ admissible pair.

One should be just a bit careful about one thing. If $d$ reaches extreme values namely $d \approx \l^2$, then 
the above representation may yield frequencies $|\xi| << \l$ and then we cannot apply our reasoning anymore.
On the other hand, if $d \approx \l^2$, then $||f||_{X^{0,\q}_{\l,d}} \approx \l^{-1} ||f||_{L^{2}_{\l,d}}$
and a straightforward computation gives all the claims listed above.

We end up the section with a linear estimate explaining how we $X^{0,\pm \q,1}$ spaces fit into the linear equation.

\begin{p1} \label{xsb5}
If $f \in X^{0,-\q,1}_\l$ then the solution of \eqref{I} $u \in X^{0,\q,1}_\l$.
\end{p1}

\begin{proof} We use $\hat{u}=(\tau-\xi^2)^{-1} \hat{f}$ to claim $u \in X^{0,\q,1}_\l$. In this process we lose track of the initial
data, but then we use $X^{0,\q,1}_\l \subset C_t L^{2}_x$ to correct the initial data as before. 

\end{proof}

\section{Definition of spaces} \label{DF}

For fixed $\l$ and $\t$ we define the spaces $Y_{\l,\t}$ as follows. If $f$ is localized at frequency $\l$ we define the norm:

\beq \label{ynorm}
||f||_{Y_{\l,\t}}=\l^{\q}||f||_{L^{\infty}_{x_{\t}}L^{2}_{t,x_{\t}'}}
\eeq  

For functions localized at frequency $\l$, by $Y_{\l}=\sum_{\t} Y_{\l,\t}$ we mean the space with norm:

\[
||f||_{ Y_{\l}} = \inf \{{\sum_{\t} ||f^{\t}||_{Y_{\l,\t}}: f = \sum_{\t} f^{\t}} \}
\]

Notice that we do not impose any localization condition on the terms $f^{\t}$. 

For function localized at frequency $\l$ we define:

\[
||f||_{\tilde{Y}_{\l,\t}}=\l^{-\frac{n-1}{2}}||f||_{L^{2}_{x_\t}L^{\infty}_{x'_\t,t}}
\]

We collect all these norms via a suppremum:

\[
||f||_{\tilde{Y}_{\l}} = \sup_{\t} ||f||_{\tilde{Y}_{\l,\t}}
\]

For fixed $\l$ we define $Z_{\l}$ to be

\beq \label{defs}
\tilde{Z}_{\l}= C_t L^{2}_x \cap Y_{\l} \cap \tilde{Y}_\l \cap L^{\frac{2(n+2)}{n}}_{x,t}  \cap X^{0,\q,\infty}_{\l}
\eeq

\noindent
with the natural assumption that the functions in $\tilde{Z}_\l$ are assumed to be localized at frequency $\l$. Then define

\beq \label{defs1}
Z_{\l}=\{ f \in \tilde{Z}_\l: f_{\cdot, \leq d} \in Z_\l, \forall d \}
\eeq

\noindent
with the norm in $Z_\l$ being defined by $||f||_{Z_\l}=\sup_{d} ||f_{\cdot, \leq d}|_{\tilde{Z}_\l}$.

From \eqref{xsb1es1}-\eqref{xsbstr} we have that 

\beq \label{xsub}
X^{0,\q,1}_{\l} \subset Z_\l 
\eeq

The total space $Z^\frac{n}{2}$ is the space of distributions satisfying:

\beq \label{distr}
f=\sum_{\l} f_\l
\eeq

\noindent
in a distributional sense and whose norm:

\[
||f||_{Z^\frac{n}{2}}=\sum_{\l} \l^{\frac{n}{2}}||f_\l||_{Z_\l}
\]

\noindent
is finite. It is a straightforward exercise to show that $Z^{\frac{n}{2}} \subset C_t \dot{B}^{2,1}_{\frac{n}{2}}$. The basic idea is to decompose

\[
f_\l=  s_{\l}(D) f_\l + \sum_{\l' < \l} s_{\l'}(D) f_\l
\]

\noindent
and notice that for $s_{\l}(D) f_\l$ we have the right estimate while for the terms $s_{\l'}(D) f_\l$ with $\l' < \l$ (which corresponds to very high modulations) we have better estimates which are enough to ensure all later summations. 

One can also show that $Z^{\frac{n}{2}}$ is a Banach space; one cannot claim this at a dyadic level, i.e. one cannot claim that $Z_\l$ is a Banach space. This is due to the localization required and not to the choice of norms. Since the structures defining $Z_\l$ and $Z_\m$ are compatible for $\l \approx \m$, this artificial problem can be easily dealt with in the global structure $Z^{\frac{n}{2}}$.

Next, we define the space for the nonlinearity, $W^{s}$. We need the dual of $Y_{\l}$, which we denote by $\mathcal{Y}_{\l}$, whose definition is:

\[
||f||_{ \mathcal{Y}_{\l}} = \sup_{\t} ||f||_{\mathcal{Y}_{\l,\t}} \ \ \ \mbox{where} \ \ \
||f||_{\mathcal{Y}_{\l,\t}}=\l^{-\q}||f||_{L^{1}_{x_{\t}}L^{2}_{t,x_{\t}'}}
\]

We define $\tilde{W}_{\l}$ the space of functions localized at frequency $\l$ by:

\[ 
\tilde{W}_{\l}=\mathcal{Y}_{\l} + L^{\frac{2(n+2)}{n+4}} + X^{0,-\q,1}_{\l}
\]

\noindent
with the standard norm. From \eqref{dual} and \eqref{xsub} it follows, by duality, that $ \tilde{W}_{\l} \subset X^{0,-\q,\infty}_{\l}$. Then we define:

\[
W_{\l, \leq d}=\{u: u=f_{\cdot, \leq d}, f \in \tilde{W}_\l \}
\]

\noindent
with the norm $||u||_{W_{\l,\leq d}}=\inf ||f||_{\tilde{W}_\l}$ where the infimum is taken over those $f \in \tilde{W}_\l$ such that $u=f_{\cdot, \leq d}$. Then we define $W_\l$ by:

\[
W_\l=\sum_{d} W_{\l, \leq d}
\]

\noindent
with the standard norm. The total space $W^{\frac{n}{2}}$ is the space of distributions satisfying \eqref{distr} and whose norm:

\[
||f||_{W^{\frac{n}{2}}}=\sum_{\l} \l^{\frac{n}{2}} ||f_\l||_{W_\l}
\]

\noindent
is finite. If we consider the equation:

\begin{equation} \label{eq1}
\begin{cases}
\begin{aligned} 
iu_{t}-\Delta u &= f \\
u(x,0) &=g
\end{aligned}
\end{cases}
\end{equation}

\noindent
then we have the linear estimate:

\beq \label{eq2}
||u||_{Z^{\frac{n}{2}}} \les ||f||_{W^{\frac{n}{2}}} + ||g||_{\dot{B}^{2,1}_{\frac{n}{2}}}
\eeq

All the ingredients needed for this estimate were provided in the previous section. We quickly recall each case to make sure of that.
We check them at a dyadic level, since after that a standard summation argument gives us \eqref{eq2}.

If $f \in X^{0,-\q,1}_\l$, then Proposition \ref{xsb5} and \eqref{xsb1es1} - \eqref{xsbstr} gives us the desired result. No complication 
is posed by the modulation cut-offs.

If we deal with $f_{\cdot,\leq d}$ for $f \in \mathcal{Y}_\l$, then one notices that whenever we 
have recovered an $L^{\infty}_{x_\t}L^{2}_{x_\t',t}$ norm we localized both $\hat{f}$ and $\hat{u}$ in $B_\t$. This is fine, since
on $f$ we have all the norms we need, while for $u$ we need the norm in $Y_{\l}=\sum_{\t} Y_{\l,\t}$.  Notice also that once we chose
to work with $S_\t f$, then we have conservation of the norm $L^{1}_{x_\t}L^{2}_{x_\t',t}$ under modulation cut-offs, see Lemma \ref{cons}.
 We still left out the part of $\hat{u}$ and $\hat{f}$ which is supported in the set $|\tau-\xi^2| \geq \frac{|\xi^2|}{10}$. Since $ \tilde{W}_{\l} \subset X^{0,-\q,\infty}_{\l}$,
it follows that $W_\l \subset X^{0,-\q,\infty}_{\l}$, therefore we can place that part of $f$ in
$X^{0,-\q,\infty}$ and, since we deal with a finite range of modulations, we can actually place it in $X^{0,-\q,1}$.

If we deal with $f_{\cdot,\leq d}$, for $f \in L^{\frac{2(n+2)}{n+4}}$, then we make use of the results in Proposition \ref{repeat1}. 
The results in that Propositions do not cover the case when we use modulation cut-offs on $f$.
We can gain that result by noticing that if $v$ solves:

\[
iv_t-\D v=f
\]

\noindent
then $v_{\cdot, \leq d}$ solves the same equation with $f$ replaced by $f_{\cdot,\leq d}$. Now, the result in Proposition \ref{repeat1} allows us to
claim the estimates we need for $v_{\cdot, \leq d}$. 

\vspace{.1in}

We would like to end the section with a commentary about the structure of $Z_\l$, more exactly about the $Y_{\l}$ part of it. Naturally, one would expect to be able to measure $S_\t f$ in $L^{\infty}_{x_{\t}}L^{2}_{t,x_{\t}'}$ and this is indeed the case, see the linear estimates. The reason 
we use a more relaxed version of this, see the sum part, has a computational explanation. If we specify that $S_\t f$ should be the part in $f$ to be placed in $L^{\infty}_{x_{\t}}L^{2}_{t,x_{\t}'}$, one runs into rather technical arguments when trying to prove that $Z^{\frac{n}{2}}+\bar{Z}^{\frac{n}{2}}$ is an algebra.  
By not being explicit about what part of $f$ we place in $L^{\infty}_{x_{\t}}L^{2}_{t,x_{\t}'}$ we simplify significantly the argument. 

\section{Bilinear Estimates} \label{BE}

The following result contains all the bilinear estimates we need.

\begin{p1} \label{mbe}
a) If $\l \leq \m$ then 

\beq \label{b1}
||B(u_\l,v_\m)_\m||_{W^{\frac{n}{2}}} \les ||u_{\l}||_{Z_\l^\frac{n}{2}} ||v_{\m}||_{Z_{\m}^\frac{n}{2}}
\eeq

\beq \label{b2}
||B(u_\m,v_\m)_\l||_{W^{\frac{n}{2}}} \les  (\l \m^{-1})^\frac{n}{2} ||u_{\m}||_{Z_\m^\frac{n}{2}} ||v_{\m}||_{Z_{\m}^\frac{n}{2}}
\eeq

\end{p1}

Once we have this result a standard summation argument gives us \eqref{bee}.

Proving the above proposition requires a good understanding of the effect of our null form. In other words, we need
a quantitative way the express the main property of our special nonlinearity. A direct computation gives us:

\[
2 \nabla{u} \cdot \nabla{v}=(i\partial_t - \D)u \cdot v + u \cdot(i\partial_t - \D)v -(i\partial_t - \D)(u \cdot v) 
\]

Our computations involve $u_{\l,\leq d_1}, v_{\m,\leq d_2}$ and  we want to 
estimate the part of $u_{\l,\leq d_1} \cdot v_{\m,\leq d_2}$ localized at modulation $\leq d_3$. 
In \cite{be2} we showed that the contribution of the gradients is 
of order $\min{(\l \cdot \m, \max{(d_1,d_2,d_3)})}$.

If $\max{(d_1,d_2,d_3)} \geq \l \m$ then we do not do anything special; 
we just use the fact that:

\[
||\nabla{u}_{\l,\leq d_1}||_{Z_\l} \leq \l ||u_{\l,\leq d_1}||_{Z_\l} 
\]

\noindent
and the similar one for $v_{\m,\leq d_2}$. If $\max{(d_1,d_2,d_3)} \leq \l \m$ then we 
would like to capture the fact that the contribution of the gradients is of
order $\max{(d_1,d_2,d_3)}$, but our spaces do not allow us to do this directly.

We proceed as follows:

\[
\begin{split}
&\ \ \ \ 2 M_{\leq d_3 } (\nabla{u}_{\l,\leq d_1} \cdot \nabla{v}_{\m,\leq d_2}) \\
&=M_{\leq d_3} \left( (i\partial_t - \D)u_{\l,\leq d_1} \cdot v_{\m,\leq d_2} + u_{\l,\leq d_1} \cdot(i\partial_t - \D)v_{\m,\leq d_2} \right) - M_{\leq d_3}(i\partial_t - \D)(u_{\l,\leq d_1} \cdot v_{\m,\leq d_2})
\end{split}
\]

For the third term
we anticipate a result (to be derived in the next section), namely that:

\[
||(u_{\l,\leq d_1} \cdot v_{\m,\leq d_2})_{\cdot, \leq \l \m}||_{X^{0,\q,1}} \leq ||u_\l||_{Z_\l} ||v_\m||_{Z_\m} 
\]

Since $(i\partial_t - \D)X^{0,\q,1}=X^{0,-\q,1}$, the third term is fine. 

For the first two terms, the key estimate is that:

\[
||(i\partial_t - \D)u_{\l,\leq d_1}||_{X^{0,\q,1}} \leq d_1 ||u_{\l,\leq d_1}||_{X^{0,\q,\infty}} 
\]

\noindent
and the similar one for $v_{\m,\leq d_2}$. Checking this estimate is a straightforward exercise. 

These observations allow us to treat the bilinear term $B(u,v)$ as follows:

- if $\max{(d_1,d_2,d_3)} \geq \l \m$ we treat $B(u,v)$ as if it were $u \cdot v$ and then factor in the gradient effect 
to be $\l \m$

- if $\max{(d_1,d_2,d_3)} \leq \l \m$ we treat $B(u,v)$ as if it were $u \cdot v$ and then factor in the gradient effect 
to be $\max{(d_1,d_2,d_3)}$

We use this simplification in order to spare space and make computations easier to follow. With this is mind we can 
move on to the proof of the bilinear estimates.

\begin{proof}[Proof of Proposition \ref{mbe}]

a) One of the estimates in the body of this proof shows that $B(u_\l,v_\m) \in L^{2}$. We need this now as a qualitative estimate, so that we
can claim:

\[
B(u_\l,v_\m)_\m=\sum_{d_3} B(u_\l,v_\m)_{\m,d_3} 
\] 

From this we can derive:

\[
B(u_\l,v_\m)_\m= \sum_{d_1} B(u_{\l, d_1},v_{\m, \leq d_1})_{\m, \leq d_1} 
\]

\[
+ \sum_{d_2} B(u_{\l, \leq d_2},v_{\m,d_2})_{\m, \leq d_2}+ \sum_{d_3} B(u_{\l, \leq d_3},v_{\m, \leq d_3})_{\m,d_3}
\]

We estimate each of these sums one at a time. For the first sum we fixed $d_1$ and note that the volume of the support of $\hat{v}_{\l,d_1}$ is $\approx \l^n d_1$. Then we estimate: 

\[
\begin{split}
||B(u_{\l,d_1},v_{\m, \leq d_1})_{\m, \leq d_1}||_{\tilde{W}_{\m,d_1}} &\les d_1 ||u_{\l,d_1} \cdot v_{\m, \leq d_1}||_{L^{\frac{2(n+2)}{n+4}}_{x,t}} \\
&\les d_{1} ||u_{\l,d_1}||_{L^{\frac{n}{2}+1}_{x,t}} ||v_{\m, \leq d_1}||_{L^{\frac{2(n+2)}{n}}_{x,t}} \\
&\les d_{1} (\l^n d_1)^{\frac{n-2}{2(n+2)}} ||u_{\l,d_1}||_{L^{2}_{x,t}} ||v_{\m, \leq d_1}||_{L^{\frac{2(n+2)}{n+4}}_{x,t}} \\
&\les d_{1}^{\q} (\l^n d_1)^{\frac{n-2}{2(n+2)}} ||u_{\l,d_1}||_{X^{0,\q}} ||v_{\m, \leq d_1}||_{L^{\frac{2(n+2)}{n+4}}_{x,t}} \\
&\les \l^{\frac{n}{2}} (\l^{-2} d_1)^{\frac{n}{n+2}} ||u_{\l,d_1}||_{Z_\l} ||v_{\m, \leq d_1}||_{Z_\mu}
\end{split}
\]

$d_1$ runs on a dyadic range up to $\approx \l^2$, hence by performing the summation with respect to $d_1$ we obtain:

\[
||\sum_{d_1} B(u_{\l, d_1},v_{\m, \leq d_1})_{\m, \leq d_1}||_{W_\m} \les \l^{\frac{n}{2}} ||u_{\l,d_1}||_{Z_\l} ||v_{\m, \leq d_1}||_{Z_\mu}
\]

To estimate the second term, we fix $d_2$ and compute: 

\[
\begin{split}
||B(u_{\l, \leq d_2},v_{\m, d_2,})_{\m, \leq d_2}||_{\tilde{W}_{\m,\leq d_2}} &\les \min{(d_{2}, \l \m)} ||u_{\l, \leq d_2} \cdot v_{\m, d_2}||_{L^{1}_{x_{\t}}L^{2}_{x_{\t}',t}} \\
&\les \min{(d_{2}, \l \m)} ||u_{\l, \leq d_2}||_{L^{2}_{x_{\t}}L^{\infty}_{x_{\t}',t}} ||v_{\m, d_2}||_{L^2} \\
&\les \l^{\frac{n-1}{2}} d_{2}^{-\q} \min{(d_2, \l \m)} ||u_{\l}||_{Z_\l} ||v_{\m}||_{Z_{\m}}
\end{split}
\]

Passing to the norm in $Y_{\l,\t}$ we obtain:

\[
||B(u_{\l, \leq d_2},v_{\m, d_2})_{\m, \leq d_2}||_{Y_{\l,\t}} \les \l^{\frac{n}{2}} \min{(\l^{-1} \m^{-1} d_2, \l \m d_2^{-1})}^{\q} ||u_{\l}||_{Z_\l} ||v_{\m}||_{Z_{\m}}
\]

Then we perform the summation with respect to $d_2$ by splitting the range in $d_2 \leq \l \m$ and $d_2 \geq \l \m$ and obtain:

\[
||\sum_{d_2} B(u_{\l, \leq d_2},v_{\m, d_2})_{\m, \leq d_2}||_{Y_{\l,\t}} \les \l^{\frac{n}{2}} ||u_{\l}||_{Z_\l} ||v_{\m}||_{Z_{\m}}
\]

For the last sum we fix $d_3$ and estimate

\[
\begin{split}
||B(u_{\l, \leq d_3},v^{\t}_{\m, \leq d_3,})_{\m,d_3}||_{L^{2}} &\les \min{(d_{3}, \l \m)} ||u_{\l, \leq d_3}||_{L^{2}_{x_{\t}}L^{\infty}_{x_{\t}',t}} ||v^{\t}_{\m, \leq d_3}||_{L^{\infty}_{x_{\t}}L^{2}_{x_{\t}',t}} \\
&\les \l^{\frac{n-1}{2}} \m^{-\q} \min{(d_3, \l \m)} ||u_{\l}||_{Z_\l} ||v^{\t}_{\m}||_{Z_{\m}}
\end{split}
\]

Summing up with respect to $\t$ and then passing to $X^{0,-\q}$ gives us:

\[
||B(u_{\l, \leq d_3},v_{\m, \leq d_3})||_{X_{\m,d_3}^{0,-\q}} \les \l^{\frac{n}{2}} \min{(\l^{-1} \m^{-1} d_3, \l \m d_3^{-1})}^{\q} ||u_{\l}||_{Z_\l} ||v_{\m}||_{Z_{\m}}
\]

We perform the summation with respect to $d_3$ in a similar fashion to the one we performed above with respect to $d_2$ and claim

\[
||\sum_{d_3} B(u_{\l, \leq d_3},v_{\m, \leq d_3})_{\m, d_3}||_{X_{\m}^{0,-\q,1}} \les \l^{\frac{n}{2}}  ||u_{\l}||_{Z_\l} ||v_{\m}||_{Z_{\m}}
\]

Bringing all the three estimates together gives us the claim in \eqref{b1}. 

b) The argument is completely similar to the one in part a). The only part that we need to argue with is that:

\[
B(u_\m,v_\m)= \sum_{\l} B(u_\m,v_\m)_\l
\]

The last estimate in part a) gives us an $L^2$ bound for $B(u_\m,v_\m)$ and this justifies the above representation. 

\end{proof}

\section{Algebra Properties} \label{AP}

The main ingredient for the algebra property is the following result. 

\begin{p1} If $\l \leq \m$ then

\beq \label{alg1}
||u_\l \cdot v_\m||_{Z_\m^{\frac{n}{2}}} \les ||u_\l||_{Z_\l^{\frac{n}{2}}} ||v_\m||_{Z_{\m}^{\frac{n}{2}}}
\eeq

\beq \label{alg2}
||u_\m \cdot v_\m||_{Z_\l^{\frac{n}{2}}} \les (\l \m^{-1})^{\frac{n}{2}} ||u_\l||_{Z_\l^{\frac{n}{2}}} ||v_\m||_{Z_{\m}^{\frac{n}{2}}}
\eeq

\end{p1}

\begin{proof}

We start the proof with showing that the part of $u_\l v_\m$ supported at modulation $\les \l \m$ can be placed in $X^{0,\q,1}$. Recall that we have made this claim in the previous section and used it there. We start by:

\[
\begin{split}
||u_{\l} \cdot v_{\m}^{\t}||_{L^{2}} &\les ||u_{\l}||_{L^{2}_{x_\t}L^{\infty}_{x_\t',t}} ||v_{\m}^{\t}||_{L^{\infty}_{x_\t}L^{2}_{x_\t',t}} \\
& \les \l^{\frac{n-1}{2}} \m^{-\q} ||u_\l||_{Z_\l} ||v_\m^{\t}||_{Y_{\l,\t}} 
\end{split}
\]

Using the definition of $Y_\l$ we conclude with:

\[
||u_{\l} \cdot v_{\m}||_{L^{2}} \les \l^{\frac{n}{2}} (\l \m)^{-\q} ||u_\l||_{Z_\l} ||v_\m||_{Z_\m} 
\]

Therefore we can conclude with:

\[
||(u_{\l} \cdot v_{\m})_{\cdot, \leq  2^{5} \l \m}||_{X^{0,\q,1}} \les \l^{\frac{n}{2}} ||u_\l||_{Z_\l} ||v_\m||_{Z_\m} 
\]

Note. The argument above is valid if we replace $u_\l$ with $u_{\l, \leq d_1}$ and $v_{\m}$ with $v_{\m,\leq d_2}$, for any $d_1, d_2$.  

Therefore we need to show only that $(u_{\l} \cdot v_{\m})_{\cdot, \geq  2^{5} \l \m} \in Z_{\m}$. First we show that $(u_{\l} \cdot v_{\m})_{\cdot, \geq  2^{5} \l \m} \in X^{0,\q,\infty}$. Since we localize at modulation higher than $2^{5} \l \m$, the low frequency cannot change 
the modulation of the high frequency too much, namely if $d \geq 2^{5} \l \m$, then:

\[
(u_{\l} \cdot v_{\m})_{\cdot, d}= ( u_\l (v_{\m, 2^{-1}d} + v_{\m,d}+ v_{\m, 2d}))_{\cdot, d}
\]

Using the $L^{\infty}$ bound on $u_\l$ and the $X^{0,\q,\infty}$ bound for the high frequency terms we obtain:

 \[
||(u_{\l} \cdot v_{\m})_{\cdot, d}||_{X^{0,\q}} \les  \l^{\frac{n}{2}} ||u_\l||_{Z_\l} ||v_\m||_{Z_\m}
\]

\noindent
and this justifies the $X^{0,\q,\infty}$ structure for the product. For all the other structures we proceed as follows.
As before if $d \geq 2^{5} \l \m$ we have:

\[
(u_{\l} \cdot v_{\m})_{\cdot, \leq d}=(u_{\l} \cdot v_{\m, \leq 2d})_{\cdot, \leq d}=u_{\l} \cdot v_{\m, \leq 2d} - (1-M_{\leq d})(u_{\l} \cdot v_{\m, \leq 2d}) 
\]

All the space-time norms are stable under multiplication by $L^{\infty}_{x,t}$ (recall $Z_{\l} \subset \l^{\frac{n}{2}} L^{\infty}$), therefore:

\[
||u_{\l} \cdot v_{\m, \leq 2d}||_{X} \les \l^{\frac{n}{2}} ||u_\l||_{Z_\l} ||v_{\m, \leq 2d}||_{X}
\]

\noindent
for any $||\cdot||_{X}$ space time norm in the definition of $Z_{\m}$ (basically all but $X^{0,\q,\infty}$). $(1-M_{\leq d})(u_{\l} \cdot v_{\m, \leq 2d})$ is supported at modulation $\approx d$, hence we have the $X^{0,\q,1}$ estimate for that term which comes from the $X^{0,\q,\infty}$
estimate we have derived before and the finite range of modulations. Hence we can conclude with:

\[
||u_{\l} \cdot v_{\m}||_{Z_\l} \les \l^{\frac{n}{2}} ||u_\l||_{Z_\l} ||v_{\m}||_{Z_\m}
\]

The proof of \eqref{alg2} follows the same lines. As before we argue that:

\[
u_\m v_\m=\sum_{\l} (u_\m v_\m)_\l
\]

\noindent
by invoking an $L^2$ estimate on the product. 
\end{proof}

Once we have this result we can derive that $Z^{\frac{n}{2}}+\bar{Z}^{\frac{n}{2}}$ is an algebra. From \eqref{alg1} and \eqref{alg2} we can conclude (via the standard summation argument) that $Z^{\frac{n}{2}} \cdot Z^{\frac{n}{2}} \subset Z^{\frac{n}{2}}$. By conjugation we obtain that $\bar{Z}^{\frac{n}{2}} \cdot \bar{Z}^{\frac{n}{2}} \subset \bar{Z}^{\frac{n}{2}}$. We are left with showing that $Z^{\frac{n}{2}} \cdot \bar{Z}^{\frac{n}{2}} \subset Z^{\frac{n}{2}} +\bar{Z}^{\frac{n}{2}}$. 

The idea is that if $\l \leq \m$, then $Z^{\frac{n}{2}}_\l \cdot \bar{Z}^{\frac{n}{2}}_\m \subset \bar{Z}^{\frac{n}{2}}_{\m}$ and if $\l \geq \m$, then $Z^{\frac{n}{2}}_\l \cdot \bar{Z}^{\frac{n}{2}}_\m \subset Z^{\frac{n}{2}}_{\l}$. To prove this, one notices again that the space-time parts of the norm do not see the difference coming from the conjugation. 
It is only the $\bar{X}^{{\frac{n}{2}},\q,\infty}$ structure which is significantly different than $X^{{\frac{n}{2}},\q,\infty}$. An argument similar to the the one in \eqref{alg1} leads to the desired result. A basic outline goes as follows:

- assume $\l \leq \m$

- $u_\l \cdot v_{\m, \leq 2^{5} \l \m}$ is placed in $L^{2}$, and then in $\bar{X}^{0,\q,1}$ 

- $u_\l \cdot v_{\m, \geq 2^{5} \l \m}$ is placed in $\bar{X}^{0,\q,\infty}$ on behalf of the $\bar{X}^{0,\q,\infty}$ structure of the high frequency, $v_{\m}$, and the modulation localization

- then we recover all the space-time norms as before.

\vspace{.1in}

Next, we turn our attention to the proof of claim \eqref{M}. 

\begin{p1} If $\l \leq \m$ then

\beq \label{alg3}
||u_\l \cdot v_\m||_{W_\m^{{\frac{n}{2}}}} \les  ||u_\l||_{Z_\l^{\frac{n}{2}}+\bar{Z}^{\frac{n}{2}}_\l} ||v_\m||_{W_{\m}^{\frac{n}{2}}}
\eeq

\beq \label{alg4}
||u_\m \cdot v_\m||_{W_\l^{\frac{n}{2}}} \les (\l \m^{-1})^{\frac{n}{2}}  ||u_\m||_{Z_\l^{\frac{n}{2}}+\bar{Z}^{\frac{n}{2}}_\l} ||v_\m||_{W_{\m}^{\frac{n}{2}}}
\eeq

If $\l \geq \m$ then

\beq \label{alg5}
||u_\l \cdot v_\m||_{W_\l^{\frac{n}{2}}} \les ||u_\l||_{Z_\l^{\frac{n}{2}}+\bar{Z}^{\frac{n}{2}}_\l} ||v_\m||_{W_{\m}^{\frac{n}{2}}}
\eeq

\end{p1}

Then a standard summation argument gives the claim in \eqref{M}.

\begin{proof}

We start with \eqref{alg3}. We can easily prove the estimate $(Z_\l+\bar{Z}_\l) \cdot X^{0,-\q,\infty}_{\m, \leq 2^{5}\l \m } \rightarrow \mathcal{Y}_{\m}$; notice that this is dual to $(Z_\l+\bar{Z}_\l) \cdot Y_\m \rightarrow X^{0,\q,1}_{\m, \leq 2^{5}\l \m }$ and we have proved this in the algebra properties. 

Therefore, what is left to be proved, is that $(Z_\l+\bar{Z}_\l)$ conserves the $W_{\m, \leq d}$ structures (for $d \geq 2^{5} \l \m$). The argument follows the same steps as in the algebra properties for the similar situation:

- we notice that $u_\l$ does not modify essentially the modulation localization of $v_{\m, d}$ for $d \geq 2^5 \l \m$

- the $X_\m^{0,-\q,\infty}$ structure is conserved (for modulations $\geq 2^5 \l \m$)  

- finally if $d \geq 2^5 \l \m$ and $v \in W_{\m, \leq d}$ then $v=f_{\cdot, \leq d}$ for $f \in \tilde{W}_\m$. For $u_\l \in Z_\l+\bar{Z}_\l$ we have that $u_\l f \in \tilde{W}_\m$ (by losing a factor of $\l^{\frac{n}{2}}$ in the estimate), hence $(u_\l f)_{\cdot, \leq d} \in \tilde{W}_{\m, \leq d}$. A direct computation shows that:

\[
M_{\leq 2^{-3}d} ((u_\l f)_{\cdot, \leq d}-u_\l f_{\cdot, \leq d})=0
\]

\noindent
therefore $u_\l f_{\cdot, \leq d}=(u_\l f)_{\cdot, \leq d}-M_{\geq 2^{-2}d} ((u_\l f)_{\cdot, \leq d}-u_\l f_{\cdot, \leq d})$. We have estimated already the first term, while the rest can be placed in $X^{0,-\q,1}$ since we established the $X^{0,-\q,\infty}$ conservation and we deal with a finite
range of modulations. 

 \eqref{alg4} is obtained in a similar way and the justification for:
 
 \[
 u_\m v_\m=\sum_\l (u_\m v_\m)_\l
 \]
 
 \noindent
 comes from the fact that in this case we obtain an $L^1_{x_\t} L^2_{x_{\t'},t}$ estimate for $u_\m v_\m$ for all $\t$ which can be easily converted into an $L^2$ estimate
 due to the frequency localization. 
 	
For \eqref{alg5} we have several cases to consider. 

Case 1. $v_\m \in X^{0,-\q,1}$. We exploit this via the obvious estimate:

\[
||v_\m||_{X^{0,\q,1}} \les \m^2 ||v_\m||_{X^{0,-\q,1}}
\]

We continue with:

\[
\begin{split}
||u_{\l} \cdot v_{\m}||_{L^{\frac{2(n+2)}{n+4}}_{x,t}} &\les ||u_\l||_{L^{\frac{2(n+2)}{n}}_{x,t}} ||v_{\m}||_{L^{\frac{n}{2}+1}_{x,t}}  \\
&\les \m^{\frac{n}{2}-2} ||u_l||_{Z_\l+\bar{Z}_\l} ||v_{\m}||_{X^{0,\q,1}} \\
&\les \m^{\frac{n}{2}} ||u_l||_{Z_\l+\bar{Z}_\l} ||v_{\m}||_{X^{0,-\q,1}}
\end{split}
\]

Case 2. $v_\m \in \mathcal{Y}_\m$. 

 Using a Sobolev embedding for the $v_\m$, we have:

\[
\begin{split}
||u_{\l} \cdot v_{\m}||_{L^{\frac{2(n+2)}{n+4}}_{x,t}} &\les ||u_\l||_{L^{\frac{2(n+2)}{n}}_{x,t}} ||v_{\m}||_{L^{\frac{n}{2}+1}_{x,t}}  \\
&\les \m^{\frac{n}{2}} ||u_\l||_{Z_\l+\bar{Z}_\l} ||v_{\m}||_{\mathcal{Y}_\m}
\end{split}
\]

Case 3. $v_\m \in L^{\frac{2(n+2)}{n+4}}$. Using a Sobolev embedding for the $v_\m$, we have:

\[
\begin{split}
||u_{\l} \cdot v_{\m}||_{L^{\frac{2(n+2)}{n+4}}_{x,t}} &\les ||u_\l||_{L^{\frac{2(n+2)}{n}}_{x,t}} ||v_{\m}||_{L^{\frac{n}{2}+1}_{x,t}}  \\
&\les \m^{\frac{n}{2}} ||u_\l||_{Z_\l+\bar{Z}_\l} ||v_{\m}||_{L^{\frac{2(n+2)}{n+4}}}
\end{split}
\]

In Case 2 and 3 we should also be aware that $v_\m$ can be of type $f_{\cdot, \leq d}$ with $f \in \mathcal{Y}_\m$ or $L^{\frac{2(n+2)}{n+4}}$. 
This can be easily dealt with by noticing that in our Sobolev embeddings we pass between two spaces which contain $L^{2}$ in between.
Therefore we can do it in two steps: use a first embedding in $L^{2}$ (where the modulation cut-off is harmless) and then pass to the $L^{\frac{n}{2}+1}_{x,t}$ norm.  

\end{proof}

\section{Smooth solutions}

In this section we prove the fact that if, in addition, our initial data $u_0 \in \dot{H}^{s}$ for $s > \frac{n}{2}$ then the solution 
is smoother. For this we define the norm in $Z^s$ by:

\[
||f||_{Z^s}^2=\sum_{\l} \l^{2s}||f_\l||^2_{Z_\l}
\]

\noindent
and similarly $W^s$. Notice that the main difference between $Z^{\frac{n}{2}}$ and $Z^s$ is the fashion in which we sum the norms of the dyadic pieces. In $Z^{\frac{n}{2}}$ we do it as in $\dot{B}_{\frac{n}{2}}^{2,1}$, i.e. in $l^1$, while in $Z^s$  we do it as in $\dot{H}^s$, i.e. in $l^{2}$. 

All these estimates \eqref{sm0}-\eqref{sm3} are direct consequences of the estimates we already proved. On dyadic pieces we already
have the estimates from previous sections and then a standard summation argument on dyadic pieces gives us the 
desired result. Note that the inequality $s > \frac{n}{2}$ gives us enough room to obtain an $l^1$ estimate for the 
norms on dyadic pieces from the $l^2$ one. 

\section{Appendix}

We end our paper with the proof of 

\[
||u||_{L^{2}_{x_{\t'}}L^{\infty}_{t,\tilde{x}'_{\t'}}} \les \l^{\frac{n-1}{2}} ||f||_{L^{\frac{2(n+2)}{n+4}}} 
\]

\noindent
where $u$ is the solution of \eqref{I}. As always we can assume $\l=1$ by rescaling. 
Our argument does one simple thing: it runs the one used in the proof Christ-Kiselev Lemma, see \cite{ck} and checks that
it works for our problem. This is why we use, essentially, the same notations. For any unproved claim, the reader may consult \cite{ck}.
Our solution is given by the Duhamel formula:

\[
u(x,t)= \int_{0}^{t} e^{-i(t-s)\D} f(x,s) ds=\int_{0}^{\infty} 1_{s \leq t} e^{-i(t-s)\D} f(x,s) ds=T(1_{s \leq t} f)(x,t)
\]

\noindent
where $T$ is the non retarded operator:

\[
T(f)(x,t)= \int_{0}^{\infty} e^{-i(t-s)\D} f(x,s) ds= e^{-it\D} \int_{0}^{\infty} e^{is\D} f(x,s) ds
\]

The operator $T$ has the desired property:

\[
||e^{-it\D} \int_{0}^{\infty} e^{is\D} f(x,s) ds||_{L^{2}_{x_\t}L^{\infty}_{x_\t',t}} \les ||\int_{0}^{\infty} e^{is\D} f(x,s)||_{L^2_{x}} \les ||f||_{L^{\frac{2(n+2)}{n+4}}}
\]

The first estimate is \eqref{se}, while the second is the standard dual Strichartz estimate, see \cite{gv}. Our task is to prove that the estimate still holds even when retarded. 

We consider a filtration $[0,t_k]_{l=1,..,l}$ of $\R^{+}$, i.e. an increasing finite sequence $(t_{k})_{n= 1,..,l}$, and construct the maximal operator:

\[
T^{*}(f)(x,t)=\max_{k} |T(1_{s \leq t_k}f)(x,t)| 
\]

We will show that $T^{*}$ is bounded from $L^{\frac{2(n+2)}{n+4}}$ to $L^{2}_{x_\t}L^{\infty}_{x_\t',t}$ and notice that since the estimates are independent of the filtration. As a consequence we can pass (via a limiting argument) to contable filtrations $[0,t_k]_{k \geq 1}$ of $\R^{+}$ and then one can pass to continuum filtrations (via a limiting argument) of type $[0,k]$ with $k \in \R$ and claim our result. 

We define $N: \R^n \times \R \rightarrow \{1,..,l\}$ by $N(t,x)=k$ if $\max_{k'} |T(1_{s \leq t_k'}f)(x,t)|=|T(1_{s \leq t_k}f)(x,t)|$; it turns out that $N$ is measurable. If we introduce $T^{N}f(x,t)=T(\chi_{[0,t_{N(x,t)}]}f)(x,t)$, then proving our result amounts to proving that 

\[
||T^{N}(f)||_{L^{2}_{x_\t}L^{\infty}_{x_\t',t}} \les ||f||_{L^{\frac{2(n+2)}{n+4}}}
\]

Without restricting the generality, we can assume $||f||_{L^{\frac{2(n+2)}{n+4}}}=1$. We define the probability measure on $\R$ by $\l(S)=\int_{S} |f|^{\frac{2(n+2)}{n+4}} dx dt$. 

There exist a collection $B^m_j$ of measurable subsets of $\R^{+}$ indexed by $m \in \N$ and $1 \leq j \leq 2^{m}$ such that:

- for each $m$, $(B^m_{j})_{j}$ is a partition of $\R^{+}$ into disjoint sets

- each $B^m_j$ is a union of precisely two sets $B^{m+1}_{j_1}$ and $B^{m+1}_{j_2}$

- $\l(B^m_j)=2^{-m}$ for all $m,j$

- each interval $[0,t_k]$ can be decomposed, modulo $\l$-null sets, as an empty, finite or countable infinite union:

\[
[0,t_k]=\cup_{i \geq 1} B^{m_i}_{j_i} \ \ \ \mbox{with} \ m_1 < m_2 < m_3 < ...
\] 

This decomposition may not be unique, but we choose and work with a fixed one. One can easily construct these sets by hand in our case, but \cite{ck} contains a more general argument for the existence of this collection. 

We define $A_{k}=\{(x,t): N(x,t)=k \}$. These are pairwise disjoint sets of $\R^n \times \R$. Then we define $\mathcal{R}$ to be the set of all pairs $(m,j,k)$ such that $B^m_j$ appears in the decomposition of $[0,t_k]$. Define:

\[
D^{m}_{j}=\cup_{(m,j,k) \in \mathcal{R}} A_k
\]

A straightforward argument shows that, for fixed $m$, the sets $D^m_j$ are disjoint. Then we define $f_j^m=\chi_{B^m_j}f$ and notice that:

\[
\chi_{[0,t_k]} f = \sum_{(m,j,k) \in \mathcal{R}} f^m_j
\]

Then we continue with:

\[
\begin{split}
T^{*}f &=\sum_{k} \chi_{A_k} T( \chi_{[0,t_k]} f )= \sum_{k} \sum_{(m,j,k) \in \mathcal{R}} \chi_{A_k} T(f^m_j) \\
&= \sum_{m} \sum_{j} \chi_{D^m_j} T(f^m_j) 
\end{split}
\]

We fix $m$ and since $D_j^m$ are disjoint

\[
\begin{split}
||\sum_{j} \chi_{D^m_j} T(f^m_j)||^2_{L^{2}_{x}L^{\infty}_{x',t}} & \les \sum_j || T f^m_j||^2_{L^{2}_{x}L^{\infty}_{x',t}} \les \sum_j ||f^m_j||^2_{L^{\frac{2(n+2)}{n+4}}} \\
&\les \sum_j 2^{-m \frac{n+4}{n+2}}  \les 2^{m(1-\frac{n+4}{n+2})}=2^{-\frac{2m}{n+2}} 
\end{split}
\]

Then we take a square root and perform the summation with respect to $m$ and obtain our claim.

\end{document}